\documentclass[12pt]{amsart}
\usepackage{amsthm,amssymb}
\usepackage{epsfig}

\newtheorem{thm}{Theorem}[section]
\newtheorem{lm}[thm]{Lemma}

\newtheorem{example}[thm]{Example}
\newtheorem{df}[thm]{Definition}

\newtheorem{assump}[thm]{Assumption}
\newtheorem{rem}[thm]{Remark}
\newtheorem{prob}[thm]{Problem}
\input epsf

\def \K {K\left(\begin{smallmatrix} p & q & r\\ a & b & c
\end{smallmatrix}\right)}

\begin{document}

\title[A New Way to Tabulate Knots]
{A New Way to Tabulate Knots}
\author[Lisa Hern\'{a}ndez and Xiao-Song Lin]{Lisa Hern\'{a}ndez and
Xiao-Song Lin}
\address{Department of Mathematics, University of California,
Riverside, CA 92521}
\email{lisah@math.ucr.edu and xl@math.ucr.edu}
\thanks{The authors are partially supported by NSF}
\begin{abstract}{We introduce a new way to tabulate knots by representing knot
diagrams using a pair of planar trees.
This pair of trees have their edges labeled by integers, they have no valence
2 vertices, and they have the same
number of valence 1 vertices. The number of valence 1 vertices of the trees is
called the {\it girth} of the knot diagram.
The classification problem of knots admitting girth 2 and 3 diagrams is
studied.
The planar tree pair representations of girth $\leq 3$
for knot diagrams in Rolfsen's table
are given.}
\end{abstract}
\maketitle

\section{Introduction}

The first successful tabulation of knots was done in the 1880's by a
Scottish physicist, Peter Guthrie Tait \cite{A}.  He was able to list
all the alternating knots up to ten crossings.  Since then, knots
have been primarily tabulated by number of crossings. In 1982, Morwen
Thistlethwaite used a computer program to generate a table of all
prime knots up to thirteen crossings, up to isomorphism and mirror
image. He found that there is only one knot with three crossings,
one knot with four crossings, two knots with five crossings, three
knots with six crossings, and so on, up to 9988 knots with thirteen
crossings.  All together he listed 12,965 prime knots.  This effort
has been continued and the most extensive knot table available
currently contains all 1,701,936 primes knots with the minimal
number of crossings not exceeding 16 (see \cite{HTW}). Work is
underway on further extension of this table to prime knots with 17
and 18 crossings.

Still, with this method, we can consider only finitely many knots at
each level of tabulation. In this paper, we develop a new way to
tabulate knots which allows us to study infinitely many knots at
each level of tabulation. In this tabulation, a knot is represented by
a pair of labeled planar trees. Each edge of this pair of planar
trees is labeled by an integer. Furthermore, these two trees have no
valence 2 vertices and the same number of valence 1 vertices. Fixing
such a pair of planar trees and varying the integer labels, we may get an
infinite family of knots and links. The main idea of this paper is
that it may be relatively easier to classify knots within such an
infinite family. In Section 2, we will describe how to get such a
pair of planar trees from a knot diagram. It will be clear that the
construction is similar to the construction of a Heegaard splitting
of a closed orientable 3-manifold.

Let us consider a pair of planar trees $(T,T')$ such that $T$
and $T'$ have no valence 2 vertices and the same number $g$
of valence 1 vertices. We will call this number $g$ the {\it girth} of
$(T,T')$. When $g=2$, then the number of edges of $T$
and $T'$ are both equal to 1. In this case, the knots
associated with such pairs $(T,T')$ are {\it double twist knots}. We will
denote knots in this family by $K(p,q)$, where $p,q\in\mathbb Z$ are
the integer labels of the only edge in $T$ and the only edge in
$T'$, respectively. We have $K(p,q)=K(q,p)$ and $K(\pm2,q)$ are
what we usually call {\it twist knots}. In Section 3, we will classify
double twist knots $K(p,q)$ with $p,q$ even using information obtained from
the Jones polynomial and the Conway polynomial. Note that these knots
$K(p,q)$ are special 2-bridge knots. So it is known that
they can be classified using
continuous fractions obtained from $p,q$. Nevertheless,
we present our direct
argument here for completeness, and for the reason that the same method will
be used in the next case.

The next case is $g=3$. Then the number of edges of $T$ and $T'$ can
only be 3 and
the only planar tree of 3 edges is the
Y-shaped tree. We will denote the knots associated with a pair of
Y-shaped trees by $K\left(\begin{smallmatrix} p & q &
r\\
a & b & c
\end{smallmatrix}
\right)$, where $p,q,r\in\mathbb Z$ are the edge labels of one
Y-shaped tree and $a,b,c\in\mathbb Z$ are the edge labels of the
other Y-shaped tree.

When $p,q,r,a,b,c$ are all even and positive, the diagram $\K$ presents a knot.
In this case, the Conway polynomial $\nabla_K$ can be
calculated explicitly. Using the explicit form of the
Conway polynomial, as well as a certain symmetry in the
diagram $\K$, we see that for a transposition $\tau\in S_3$,
$$\nabla_{\K}= \nabla_{K\left(\begin{smallmatrix} p & q & r\\
\tau(a) & \tau(b) & \tau(c)
\end{smallmatrix}\right)}$$
iff $$\K=K\left(\begin{smallmatrix} p & q & r\\ \tau(a) & \tau(b) &
\tau(c)
\end{smallmatrix}\right).$$

In general, the calculation of the Jones polynomial $J_K$ via the
Kauffman bracket is more manageable for $\K$ than the
Conway polynomial. In particular, the difference
$$J_{\K}-J_{K\left(\begin{smallmatrix} p & q & r\\ \tau(a) & \tau(b) & \tau(c)
\end{smallmatrix}\right)}$$
turns out to be quite simple. Thus, we will be able to show that in some
other cases, the
Jones polynomial can be used to distinguish knots in this family.
See Section 4 for these results.

Finally, we give a list of knots in the Rolfsen knot table \cite{R}
whose diagram of minimal number of crossings can be put in the form
of $K(p,q)$ or $\K$. It turns out that all knots with no more
than 7 crossings can be put into these forms, and there is only one
knot with 8 crossings, $8_{18}$, that may not be able to put into
this form. In fact, knots and links that admit planar tree pair
representations of
girth $\leq 3$ occupy a quite large portion of Rolfsen's table.

Our method is in the same spirit as Conway's method of knot
tabulation using tangles \cite{Conway}. Such methods are more
structural in their enumeration of knot diagrams. On the other hand,
the traditional method of knot tabulation using crossing
information, as invented by Dowker \cite{DowThis1, DowThis2}, has
the advantage of lending itself more easily to computer practice.
And as we mentioned before, this advantage has been fulfilled with
rather speculative success. With such a comparison in mind, it is
worthwhile to note that our method also lends itself easily to
computer enumeration of knot diagrams.

Special infinite families of knots, like torus knots, alternating
knots, and rational knots, etc., are very favorable to knot
theorists. Such families of knots all have certain rigidity that
some topological quantities are determined by their diagrammatic
descriptions. It seems to us that $\K$ may give us another such a
family of knots.
More importantly, the construction of this infinite
family of knots $\K$ can be generalized to produce systematically
other infinite families of knots, and these infinite families of
knots can be used to exhaust the entire collection of knot types.

\section{Heegaard Decomposition of a Knot}

Given a knot diagram $D_{\mathcal{K}}$ for the knot $\mathcal{K}$,
we have the checkerboard coloring of complementary regions of
$D_{\mathcal{K}}$. Recall that $D_{\mathcal{K}}$ is a generic
immersed circle in the plane with the crossing information specified
at each double point of the knot diagram. We color the complementary
regions of the knot diagram $D_{\mathcal{K}}$ by black or white such
that two regions sharing a common edge will be colored differently.
For such a checkerboard coloring, an associated planar graph
$G(D_{\mathcal{K}})$ can be defined. The set of vertices of
$G(D_{\mathcal{K}})$ will be the complementary regions of
$D_{\mathcal{K}}$ that are colored black, and there is an edge
connecting two black regions if they contain a common double point
on the diagram $D_{\mathcal{K}}$. See Figure 1.

Through out this paper, we will assume that the knot diagram $D_{\mathcal{K}}$
is {\it reduced} in the sense that $G(D_{\mathcal{K}})$ has no vertices of
valence 1.

\begin{figure}
\centerline{\epsfig{file=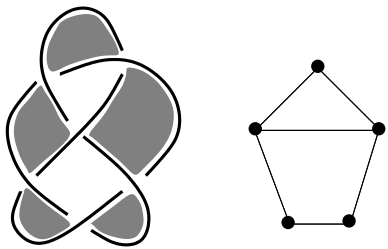,width=2in}} \caption{A knot
diagram, its checkerboard coloring, and the associated planar
graph.}
\end{figure}

Without loss of generality, we may assume that the knot diagram
$D_{\mathcal{K}}$ is contained inside of a small neighborhood of the
planar graph $G(D_{\mathcal{K}})$. Now take a spanning tree $T$ of
the graph $D_{\mathcal{K}}$, and let $U$ be a closed neighborhood of
$T$ in the plane such that

(1) $\partial U$ intersects $D_{\mathcal{K}}$ transversely;

(2) crossings of $D_{\mathcal{K}}$ in $U$ are in bijection with
edges in $T$; and

(3) the number of intersection points of $\partial U$ and
$D_{\mathcal{K}}$ is minimal among all choices of $U$ satisfying (1)
and (2) above.

Outside of $U$, we have a tree $T'$ dual to $T$, which is a spanning
tree of the dual graph $G'(D_{\mathcal{K}})$ of the planar graph
$G(D_{\mathcal{K}})$. The crossings of $D_{\mathcal{K}}$ not in $U$
are in bijection with edges in $T'$. See Figure 2.

\begin{figure}
\centerline{\epsfig{file=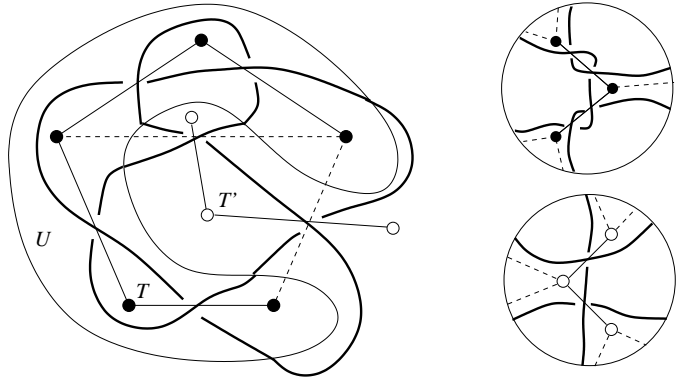,width=3.5in}} \caption{A Heegaard
decomposition of a knot diagram.}
\end{figure}

Now we can form a {\it Heegaard decomposition} of the knot diagram
$D_{\mathcal{K}}$: Redraw $D_{\mathcal{K}}$ in $U$ and outside of
$U$ in a unit disk. For the tree $T$, we omit all vertices of
valance 2 which are not connected to $\partial U$ by edges of
$G(D_{\mathcal{K}})$ not in $T$. The resulting planar tree is called
$\bar{T}$. Similarly, we can get a planar tree $\bar{T}'$ from $T'$.
Then redrawings of the diagram $D_{\mathcal{K}}$ in $U$ and outside
of $U$ both have a sequence of half twists (crossings)
on each edge of $\bar{T}$ and
$\bar{T}'$. This integer number of half twists gives a {\it label} to that
edge of $\bar{T}$ or $\bar{T}'$. See the right hand side of Figure
2, where two unit disks with $\bar{T}$ and $\bar{T}'$, as well as
the redrawings of $D_{\mathcal{K}}$ in $U$ and outside of $U$ are
shown. The original diagram $D_{\mathcal{K}}$ can be then recovered
by gluing these two unit disks together appropriately along the
boundary. In the notation specified in Section 4, this knot is
denoted by $K\left(\begin{smallmatrix} 0 & 2 & 2\\ 0 & -1 & -1
\end{smallmatrix}
\right).$

Notice that some vertices of $\bar{T}$ may be connected to $\partial
U$ by edges of $G(D_{\mathcal{K}})$ not in $T$. This is also true
for $\bar{T}'$. In the redrawings of $D_{\mathcal{K}}$ in the unit
disks, we see parts of these edges as dashed lines. For dashed lines
coming out of a vertex of $\bar{T}$ or $\bar{T}'$, the edges of
$\bar{T}$ or $\bar{T}'$ coming out of that vertex will divide them
into equivalence classes. It is clear that $\bar{T}$ and $\bar{T}'$
have the same number of equivalence classes of dashed lines. We call
this number the {\it girth} of this pair of trees
$(\bar{T},\bar{T}')$. The example in Figure 2 has the girth equal to
3. See also Figure 3 (a).

An edge in $\bar{T}$ or $\bar{T}'$ is called an {\it exterior edge}
if it has a vertex of valence 1. We can remove or add an exterior edge
with label $0$.
See Figure 3. This is how we should understand the $0$'s in the
notation like $K\left(\begin{smallmatrix} 0 & 2 & 2\\ 0 & -1 & -1
\end{smallmatrix}
\right)$. In that example, the edge number $e$ of $\bar{T}$ and $e'$
of $\bar{T}'$ are both 2, and the girth is 3. We need to add an
exterior edge to both $\bar{T}$ and $\bar{T}'$ with label 0 in order
to have the notation $K\left(\begin{smallmatrix} 0 & 2 & 2\\ 0 & -1
& -1
\end{smallmatrix}
\right)$. Thus, by adding exterior edges with label 0 if necessary,
we can make the following assumption:

\begin{assump} The girth of $(\bar{T},\bar{T}')$ is equal to the number of
valence 1 vertices of $\bar{T}$ or $\bar{T}'$, respectively.
\end{assump}

\begin{figure}
\centerline{\epsfig{file=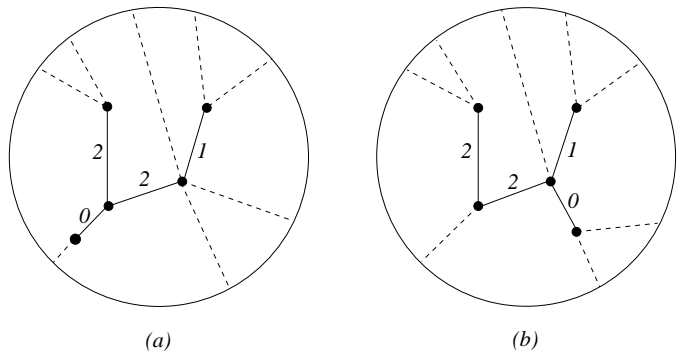,width=3.5in}} \caption{(a) The
girth in this example is 5. (b) Exterior edges with label $0$ can be
removed or added.}
\end{figure}

Finally, if the girth of the pair $(\bar{T},\bar{T}')$ is $g$, the
boundary of the unit disk ($\partial U$) is a union of arcs
$A_1,B_1,\dots,A_g,B_g$, whose interiors are all disjoint.
Furthermore, $A_1,B_1,\dots,A_g,B_g$ lie on $\partial U$ in the
given cyclic order, such that each $A_i$ contains the ends of the
dashed lines in a single equivalence class in $\bar{T}$, and each
$B_i$ contains the ends of the dashed lines in a single equivalence
class in $\bar{T}'$.

This completes our description of the essential features of Heegaard
decompositions of a knot diagram. In the next two sections, we are
going to study two cases in some detail. In the first case, we have
the pair $(\bar{T},\bar{T}')$ with the girth equal to 2 and the edge
number $e=e'=1$. In the second case, we have the pair
$(\bar{T},\bar{T}')$ with the girth equal to 3 and the edge number
$e=e'=3$.

We make the following definition.

\begin{df} The {\em girth} of a knot or link $K$ is the minimal girth
of all Heegaard decompositions of knot diagrams of $K$.
\end{df}

\section{Knots with girth 2 diagrams}

Notice the only knots with girth 2 diagrams are the {\it double
twist knots} $K(p,q)$ where there are $p$ crossings (half twists) on
one tree and $q$ on the other. See Figure 4 for the example of
$K(3,-2)$. Figure 4 also fixes our convention for positive and
negative crossings with respect to edges of the trees $T$ and $T'$.
The sign of crossings is used to determined the integer labels of
edges of $T$ and $T'$.

\begin{figure}
\centerline{\epsfig{file=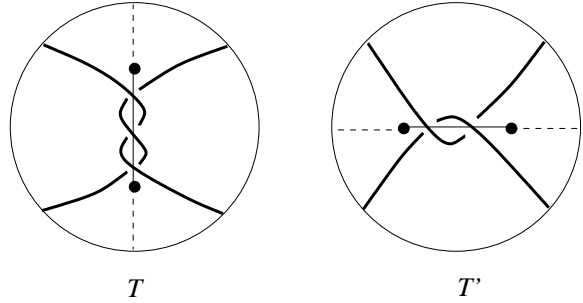,width=3in}}
\caption{The knot $K(3,-2)$.}
\end{figure}

\begin{lm} We have $K(p,q)=K(q,p)$.
\end{lm}

\begin{proof} Simply switch the trees $T$ and $T'$, we can see that the
knots $K(p,q)$ and $K(q,p)$ are isotopic.
\end{proof}

There is a special case of a {\it single
twist knot} with $p$ crossings as shown in Figure 5.
We will denote it as $K(p)$.

\begin{figure}
\centerline{\epsfig{file=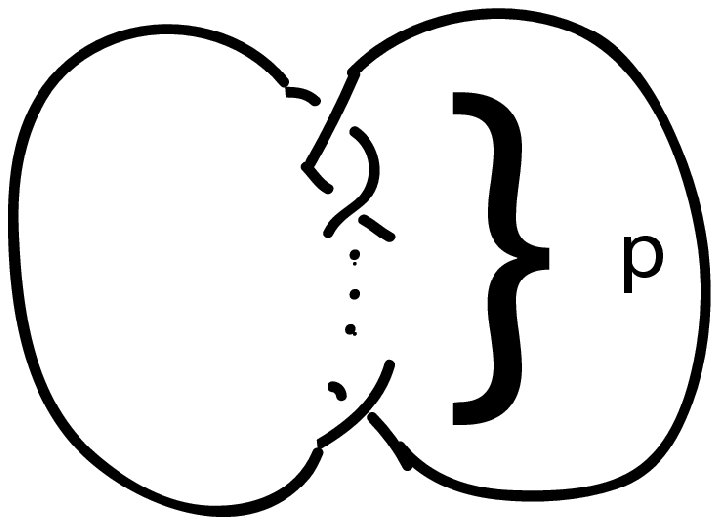,width=1in}} \vskip .3in
\centerline{\epsfig{file=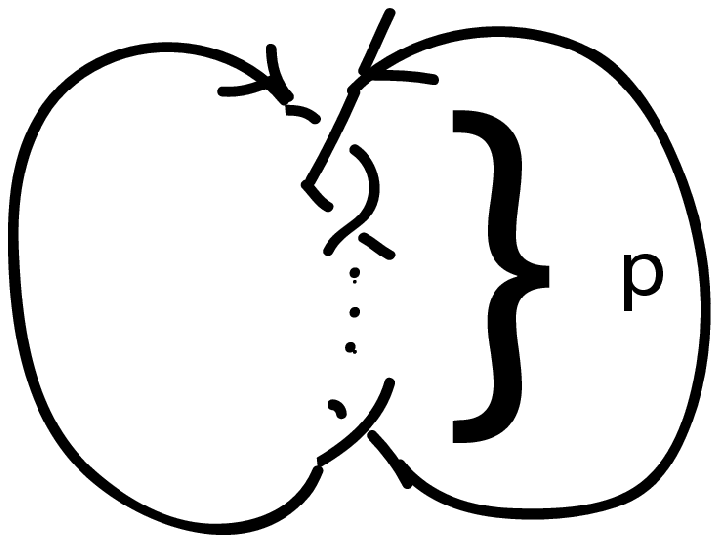,height=0.7in} \qquad \qquad
\epsfig{file=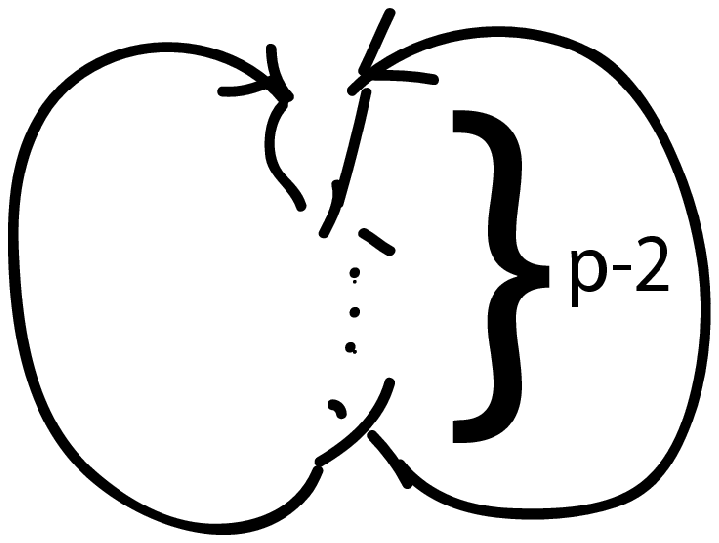,height=0.7in} \qquad \qquad
\epsfig{file=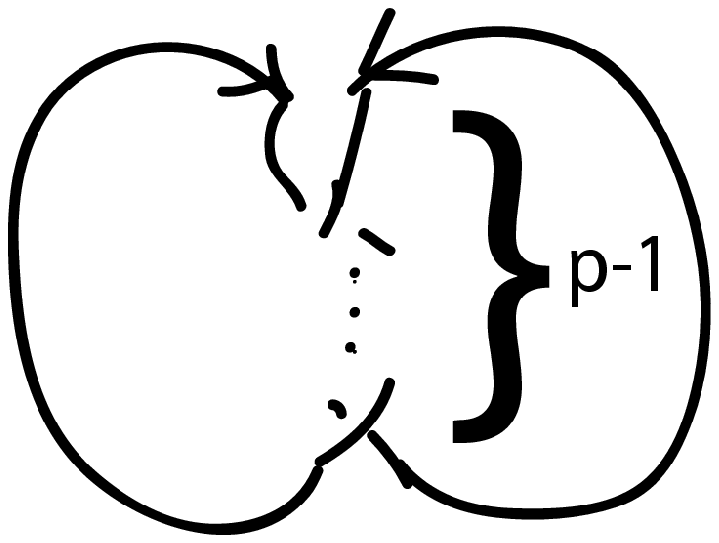,height=0.7in}}
\medskip
\centerline{\small $L_{+}$ \qquad \qquad \qquad \qquad \qquad $L_{-}$ \qquad
\qquad\qquad
\qquad \qquad $L_{0}$}
\caption{The knot $K(p)$ and its Conway skein.}
\end{figure}

\begin{lm} We have $K(p,\pm1)=K(p\mp1)$.
\end{lm}

\begin{proof} Simply look at the diagram of $K(p,\pm1)$ and we will see that
the only crossing on $T'$ can be moved to a crossing on $T$ with the sign
changed.
\end{proof}

We want to calculate the Conway polynomial $\nabla_K$ for $K=K(p)$ and
$K=K(p,q)$. Recall the Conway polynomial $\nabla_K$ is an invariant of
oriented link $K$. For $K(p)$, we want to illustrate the effect of
two different ways to orient its component(s) on its Conway polynomial.

By $K(p)$ with the same direction, we mean to orient $K(p)$ in such
a way that both strands of the chain of crossings are pointing in
the same direction. By $K(p)$ with opposite directions, we mean to
orient $K(p)$ in such a way that the strands of the chain of
crossing are pointing in different orientations. The latter case is
possible only when $p$ is even.

In the first case of $K(p)$ with the same direction, we may apply the
Conway skein relation to $L_+=K(p)$, $L_-=K(p-1)$, and $L_0=K(p-2)$
(see Figure 5) and get
$$\nabla_{K(p)} = z \nabla_{K(p-1)} + \nabla_{K(p-2)} $$
for $p\geq2$.

Now, for simplicity, denote $\nabla_{K(p)} (z)$ as $\nabla_p(z)$.

\begin{lm}\label{1-twist-1}
For $p>0$ and $K(p)$ with the same direction, we have
 $$\nabla_p(z) = i^{p-1} U_{p-1} \left( -\frac{zi}{2} \right),$$ where
$U_n$ is the n\textsuperscript{th} Chebyshev polynomial of the
second kind defined by initial values $U_0(x)=1$ and
$U_1(x)=2x$ and by the recursion $U_{n+1}(x)=2xU_n(x)-U_{n-1}(x)$, and $i=\sqrt{-1}$.
\end{lm}

If $p$ is odd, then $K(p)$ is a knot. Since the Conway
polynomial is independent of orientations for knots, it
does not matter whether $K(p)$ has the same or opposite directions.
Further, if orientation does not
matter, then the sign of $p$ does not matter. So for odd $p$, we
have $\nabla_p(z)=\nabla_{-p}(z)$. For $p$ even and $K(p)$ with
the same direction, we have $\nabla_p(z)=-\nabla_{-p}(z)$.

\begin{lm}\label{1-twist-2} For $p$ even and $K(p)$ with opposite directions, we have
$$\nabla_p (z) =\frac{p}{2}\,z \, , \mbox{\qquad if $p>0$} $$
and
$$\nabla_p (z) =-\frac{p}{2}\,z \, , \mbox{\qquad if $p<0$}.$$
\end{lm}

To compare the results of Lemma \ref{1-twist-1} and Lemma \ref{1-twist-2},
we quote the following explicit formula for the Chebyshev polynomial of
the second kind $U_p$:
$$U_p(x)=\sum_{m=0}^{[p/2]}\,\left(\begin{matrix} p+1 \\ 2m+1\end{matrix}
\right)\,x^{p-2m}(x^2-1)^m.
$$

Next, we want to calculate the Conway polynomial for $K(p,q)$.

Consider all the possibilities for values of $p$ and $q$ with all
possible orientations for double twist knots $K(p,q)$, also take into
consideration the symmetry $K(p,q)=K(q,p)$,
we see that there are the following two cases.

\begin{itemize}
\item[{\it Case 1:}] The $p$-crossings have the same direction and the
$q$-crossings have opposite directions.  Then $p$ can be either even
or odd, and $q$ must be odd.

\medskip
\centerline{\epsfig{file=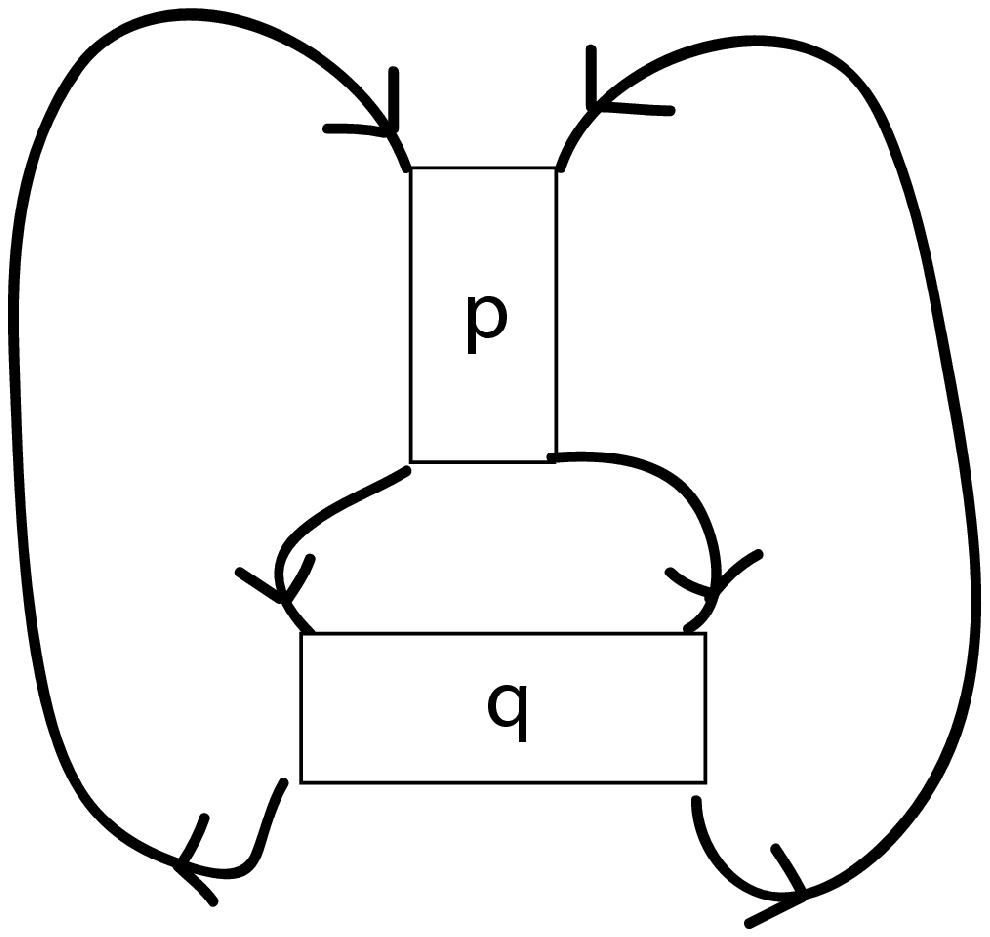,width=1.3in}}
\medskip

\item[{\it Case 2:}] Both $p$-crossings and $q$-crossings have opposite
directions. Then both $p$ and $q$ must be even. Note in this case, $K(p,q)$
must be a knot.

\medskip
\centerline{\epsfig{file=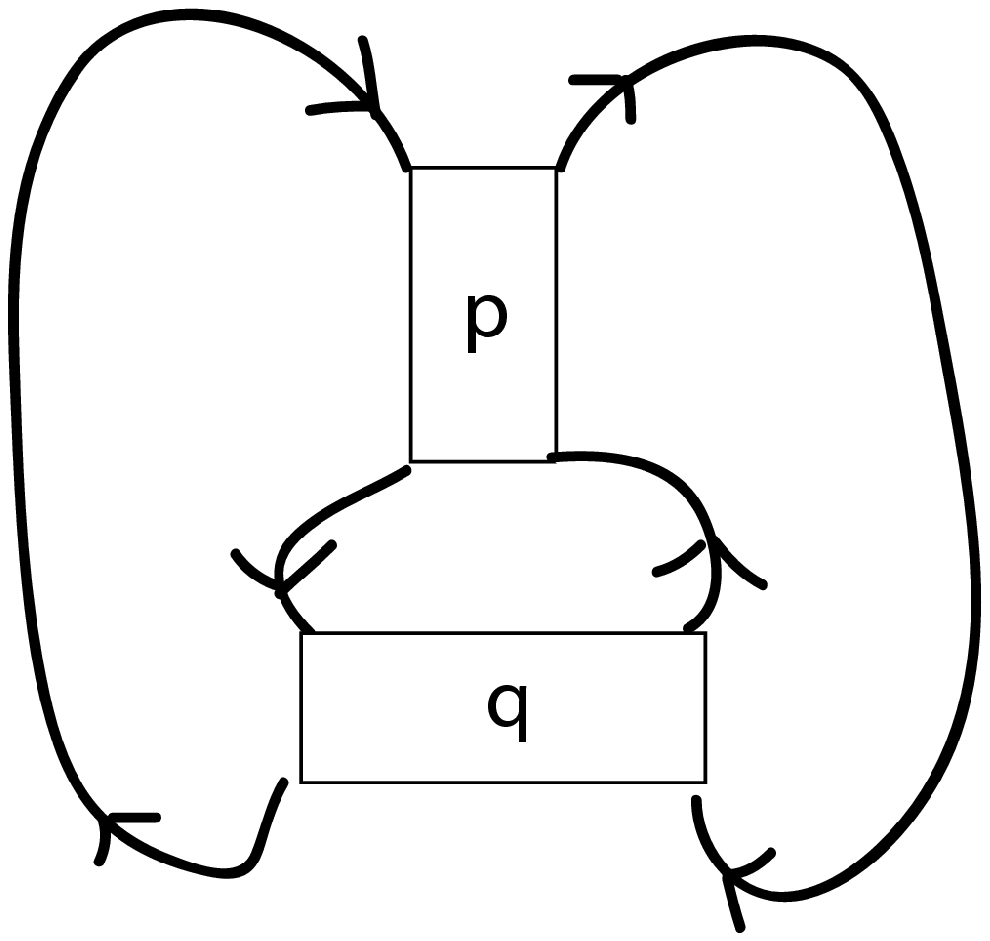,width=1.3in}}
\medskip

\end{itemize}

Since the calculation is straightforward using Conway skein relation,
we will only list the results of the calculation
in these two cases in the following lemmas.

\begin{lm}\label{conway1}
For the double twist knot $K(p,q)$ with $p$ nonzero, $q$
odd, the p-crossings having the same direction and
the q-crossings having opposite directions, the
Conway polynomial is given by
$$\nabla_{K(p,q)}(z) = \left( \frac{q-1}{2} \right) z\nabla_p(z)
+\nabla_{p+1}(z)\mbox{ \qquad if $q<0$} $$
and
$$ \nabla_{K(p,q)}(z) = \nabla_{p-1}(z) - \left( \frac{q+1}{2} \right) z
\nabla_p(z)\mbox{ \qquad if $q>0$}.$$
Here $\nabla_p(z)$ is given in Lemma \ref{1-twist-1}.
\end{lm}

\begin{lm}\label{conway2}
For the double twist knot $K(p,q)$ with $p,q$ even, the Conway
polynomial is given by
$$ \nabla_{K(p,q)} (z) = \text{\em sign$(p)$}\,\text{\em sign$(q)$}\,
\frac{pq}{4}\, z^2+1.$$
\end{lm}

In order to compute the Jones polynomial for $K(p,q)$, we define the following
Laurent polynomial of the variable $A$ (the variable in the Kauffman bracket).

\begin{df}
For $p>0$, we define
$$S_p=S_p(A)= \sum_{i=1}^p A^{2-i}(-A^3)^{p-i}.$$
And we use $S_{-p}(A) = S_p(A^{-1})$ to define $S_p$ for $p<0$.
\end{df}

Let $\langle K(p,q)\rangle$ be the Kauffman bracket of $K(p,q)$. The
following formula is easy to obtain.

\begin{lm} We have
$$\langle K(p,q)\rangle = (-A^2-A^{-2})(S_pA^{-q}+S_qA^{-p})
+S_pS_q+A^{-p-q}.$$
\end{lm}

\begin{lm}\label{jones-span} If $p,q>1$, then the lowest degree term of
$\langle K(p,q)\rangle$ is $-A^{-p-q}$, and the highest degree term
of $\langle K(p,q)\rangle$ is $(-1)^{p+q}A^{3(p+q)-4}$. Thus, the span of the
Jones polynomial of $K(p,q)$ is $p+q$.
\end{lm}

\begin{proof} Notice that the lowest degree term of $S_pS_q$ is $A^{4-p-q}$,
and the lowest degree terms of $S_pA^{-q}$ and $S_qA^{-p}$ are both
$A^{2-p-q}$. Thus the lowest degree term of
$\langle K(p,q)\rangle$ is $-A^{-p-q}$.

Similarly, notice that the highest degree term of $S_pS_q$ is $(-1)^{p+q}
A^{3(p+q)-4}$, and the highest degree terms of $S_pA^{-q}$ and $S_qA^{-p}$ are,
respectively, $(-1)^{p-1}A^{3p-2-q}$ and $(-1)^{q-1}A^{3q-2-p}$. We have
$$3(p+q)-4-3p+q=4q-4>0\quad\text{and}\quad
3(p+q)-4-3q+p=4p-4>0.$$
Therefore, the highest degree term of $\langle K(p,q)\rangle$
is $(-1)^{p+q}A^{3(p+q)-4}$.

We obtained the Jones polynomial by first orienting $K(p,q)$ and
calculating the writhe, then normalizing the Kauffman bracket by a factor
$A^w$, where the power $w$ is obtained from the writhe, and finally making the
substitution $A^4=t$. Thus the span of the Jones polynomial of $K(p,q)$ is
$$(3(p+q)-4-(-p-q))/4+1=p+q.$$
\end{proof}

\begin{example} {\em
Consider the knots $K(2,8)$ and $K(4,4)$. They have the same
Conway polynomial $1+4z^2$. But the spans of the Jones polynomial of these
two knots are, respectively, 10 and 8. Therefore $K(2,8)\neq K(4,4)$.}
\end{example}

\begin{thm}\label{pq-classification}
Suppose $p,q,a,b$ are even positive integers.
Then $K(p,q)=K(a,b)$ iff $\{p,q\}=\{a,b\}$.
\end{thm}

\begin{proof} Since $K(p,q)=K(q,p)$, $\{p,q\}=\{a,b\}$ implies $K(p,q)=K(a,b)$.

Suppose $K(p,q)=K(a,b)$. Then by Lemma \ref{conway2}, we
have $pq=ab$, and by Lemma \ref{jones-span}, we have $p+q=a+b$.
Thus, $\{p,q\}$ and $\{a,b\}$ are both equal to the pair of roots of
the quadratic equation $(x-p)(x-q)=(x-a)(x-b)=0$. Therefore
$\{p,q\}=\{a,b\}$.
\end{proof}

\begin{rem}{\rm (1) The mirror image of $K(p,q)$ is $K(-p,-q)$,
so if we do not distinguish knots and their mirror images,
by Theorem \ref{pq-classification},
knots $K(p,q)$ with $p,q$ even and $pq>0$ are classified by $\{|p|,|q|\}$.

(2) If $pq<0$, $K(p,q)$ is an alternating knot. The span of the
Jones polynomial in this case is $|p|+|q|$. So if we do not distinguish
knots and their mirror images, knots $K(p,q)$ with $p,q$ even and $pq<0$
are also classified by $\{|p|,|q|\}$.
}
\end{rem}

\section{Knots with girth 3 diagrams}

The only girth 3 tree with no valence 2 vertices is
{\epsfig{file=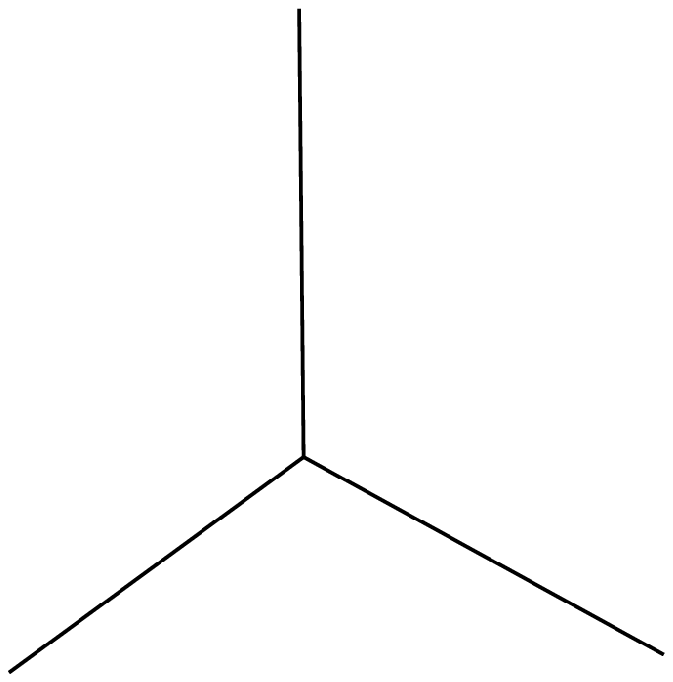,height=0.2in,width=0.2in}}. So we can
denote a knot diagram admitting a Heegaard decomposition of girth 3
as $\K$ where $p,q,r,a,b,c$ correspond to the crossings on edges of
$T$ and $T'$as depicted in Figure 6.

For example, the classical $(p,q,r)$ pretzel knots are of the form
$K\left(\begin{smallmatrix} p & q & r\\
\pm1 & \pm1 & 0\end{smallmatrix}\right)$ in our representation.

\begin{figure}
\centerline{\epsfig{file=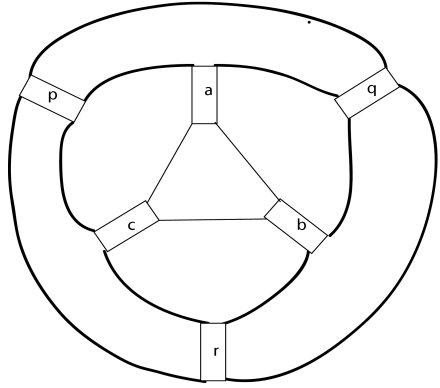,height=2in,width=2.25in}}
\caption{The knot diagram of girth 3.}
\end{figure}

The knot diagram $\K$ has some obvious symmetries. We list these symmetries
in the following lemma.

\begin{lm}
(1) $\K=K\left(\begin{smallmatrix} a & b & c\\
p & q & r
\end{smallmatrix}\right)$;

(2) $\K=K\left(\begin{smallmatrix} q & r & p\\
b & c & a\end{smallmatrix}\right)=
K\left(\begin{smallmatrix} r & p & q\\
c & a & b
\end{smallmatrix}\right)$;

(3) $\K=K\left(\begin{smallmatrix} p & r & q\\
c & a & b
\end{smallmatrix}\right)$.
\end{lm}

\begin{proof} The relation (1) corresponds to the isotopy that turns the inner ring $\{a,b,c\}$ of the knot diagram in Figure 6 to the out ring $\{p,q,r\}$
and vice versa.

For the relations (2) and (3), notice that the symmetric group of the
equilateral triangle

\centerline{\epsfig{file=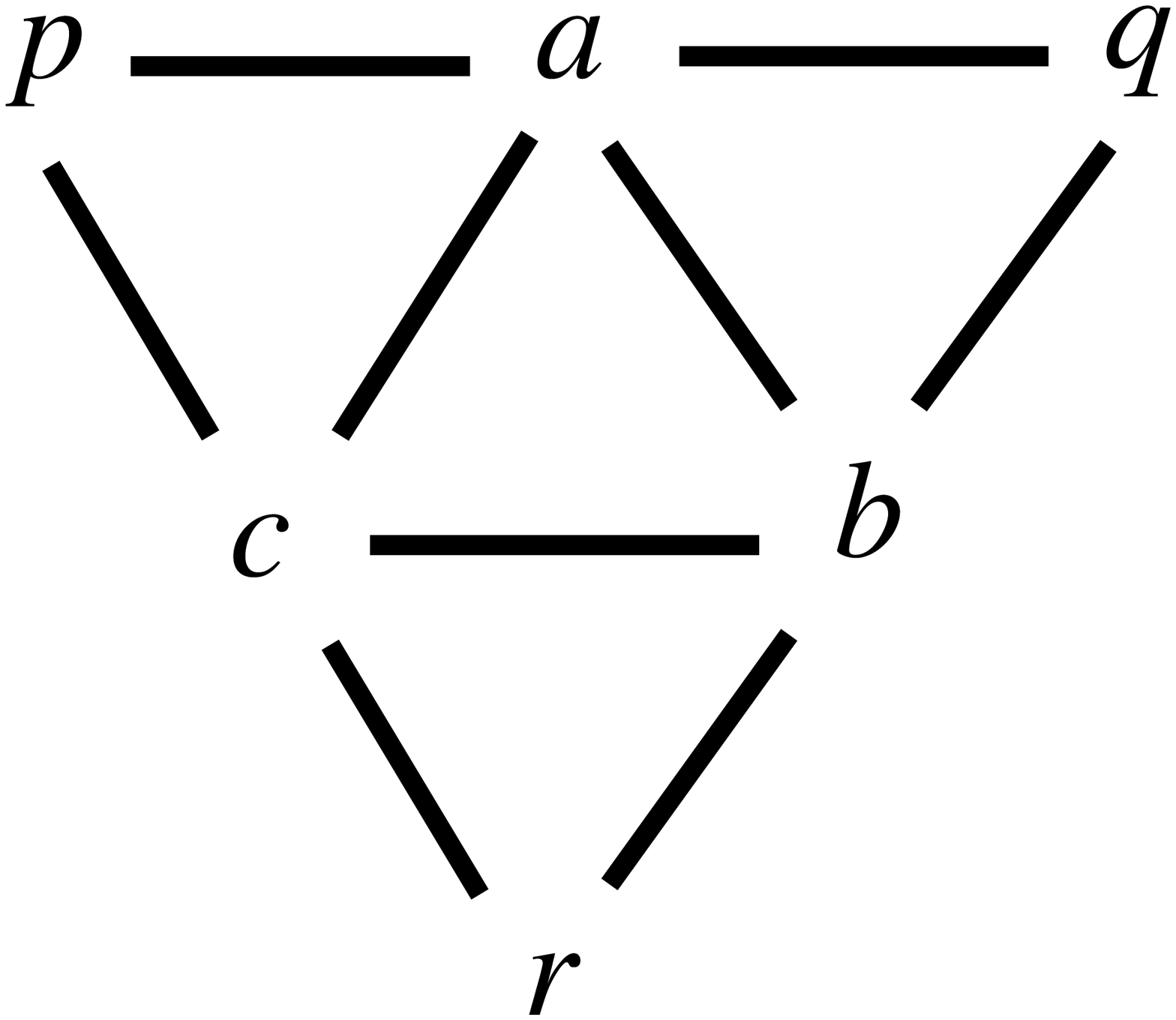,width=.7in}}

\noindent
is the dihedral group $D_3$. Every $D_3$ symmetry can be realized by an
isotopy of the knot diagram in Figure 6. Then, the relation (2) corresponds
to a rotation of $120^\circ$, and the relation (3) corresponds to a reflection.
\end{proof}

We say that $\K=K\left(\begin{smallmatrix} p' & q' & r'\\
a' & b' & c'
\end{smallmatrix}\right)$ by a $D_3$ symmetry, if $\K$ and
$K\left(\begin{smallmatrix} p' & q' & r'\\
a' & b' & c'
\end{smallmatrix}\right)$ are interchangeable by the relations
(2) and (3) in the lemma above.

Thus, the first problem in the classification of these knots $\K$
is to decide if
$$\K=K\left(\begin{smallmatrix} p & q & r\\
\tau(a) & \tau(b) & \tau(c)
\end{smallmatrix}\right)$$
for $\tau\in S_3$, the symmetric group of 3 elements.

\subsection{The Conway polynomial for $p,q,r,a,b,c$ even}

In the case when $p,q,r,a,b,c$ are all even, $\K$ is a knot. When we
orient this knot, all the crossings have opposite directions. So, as before,
the Conway polynomial has a simpler form in this case.

\begin{lm} For $p,q,r,a,b,c$ even and positive, the Conway polynomial of
$\K$ is
\begin{align}
\nabla_{K\left(\begin{smallmatrix} p & q &r\\
a & b & c
\end{smallmatrix}\right) } =&
(pq + pr + qr)(ab + ac + bc)\left(\frac{z}{2}\right)^4 \notag \\
& +(pa + pc + qa + qb + rb + rc)\left(\frac{z}{2}\right)^2 +1. \notag
\end{align}
\end{lm}

Thus, we have

\begin{equation}
\nabla_{\K}-\nabla_{K\left(\begin{smallmatrix} p & q & r\\
b & a & c
\end{smallmatrix}\right)}=(p-r)(a-b)\,\left(\frac{z}2\right)^2
\end{equation}

\begin{equation}
\nabla_{\K}-\nabla_{K\left(\begin{smallmatrix} p & q & r\\
a & c & b
\end{smallmatrix}\right)}=(p-q)(c-b)\,\left(\frac{z}2\right)^2
\end{equation}

\begin{equation}
\nabla_{\K}-\nabla_{K\left(\begin{smallmatrix} p & q & r\\
c & b & a
\end{smallmatrix}\right)}=(q-r)(a-c)\,\left(\frac{z}2\right)^2
\end{equation}

From Equations (1), (2), and (3), we can conclude the following result.

\begin{lm} Let $\tau\in S_3$ be a transposition. For $p,q,r,a,b,c$
even and positive,
$$\nabla_{\K}=\nabla_{K\left(\begin{smallmatrix} p & q & r\\
\tau(a) & \tau(b) & \tau(c)
\end{smallmatrix}\right)}$$
iff
$$\K=K\left(\begin{smallmatrix} p & q & r\\
\tau(a) & \tau(b) & \tau(c)
\end{smallmatrix}\right)$$
by a $D_3$ symmetry.
\end{lm}

\begin{proof} Suppose
$$\nabla_{\K}=\nabla_{K\left(\begin{smallmatrix} p & q & r\\
c & b & a
\end{smallmatrix}\right)}.$$
Then by Equation (3), we have either $q=r$ or $a=c$. When $q=r$, we use
Lemma 4.1 (3) to see that
$$K\left(\begin{smallmatrix} p & q & r\\
c & b & a
\end{smallmatrix}\right)=
K\left(\begin{smallmatrix} p & r & q\\
a & b & c
\end{smallmatrix}\right)=
\K.$$
When $a=c$, we certainly have
$$K\left(\begin{smallmatrix} p & q & r\\
c & b & a
\end{smallmatrix}\right)=\K.$$

Other two cases can be proved similarly.
\end{proof}

When $\tau =(132)$ or $\tau=(123)$, the Conway polynomial sometimes
fails
to distinguish $\K$ and $K\left(\begin{smallmatrix} p & q & r\\
\tau(a) & \tau(b) & \tau(c)
\end{smallmatrix}\right)$. We have the following lemma.

\begin{lm}
For $p,q,r,a,b,c$
even and positive, we have
\begin{equation}
\nabla_{\K}-\nabla_{K\left(\begin{smallmatrix} p & q & r\\
c & a & b
\end{smallmatrix}\right)}=\left|\begin{matrix}
p & q & r \\
c & a & b \\
1 & 1 & 1\end{matrix}\right|
\,\left(\frac{z}2\right)^2
\end{equation}
and
\begin{equation}
\nabla_{\K}-\nabla_{K\left(\begin{smallmatrix} p & q & r\\
b & c & a
\end{smallmatrix}\right)}=-\,\left|\begin{matrix}
p & q & r \\
a & b & c \\
1 & 1 & 1\end{matrix}\right|\left(\frac{z}2\right)^2
\end{equation}
\end{lm}

Therefore, it is possible to have
$\nabla_{\K}=\nabla_{K\left(\begin{smallmatrix} p & q & r\\
c & a & b
\end{smallmatrix}\right)}$ without having $\K=
K\left(\begin{smallmatrix} p & q & r\\
c & a & b
\end{smallmatrix}\right)$, and to have
$\nabla_{\K}=\nabla_{K\left(\begin{smallmatrix} p & q & r\\
b & c & a
\end{smallmatrix}\right)}$ without having $\K=
K\left(\begin{smallmatrix} p & q & r\\
b & c & a
\end{smallmatrix}\right)$.

\subsection{The Kauffman bracket for general $p,q,r,a,b,c$}

The computation of the Jones polynomial via the Kauffman bracket turns out to
be simpler than the Conway polynomial for $\K$.

We define first the following symmetric functions on a triplet $(p,q,r)$ as
\begin{align}
& S^1(p,q,r)=S_pA^{-q-r} +S_qA^{-p-r} +S_rA^{-p-q} \notag \\
& S^2(p,q,r)=S_pS_qA^{-r} +S_pS_rA^{-q} +S_qS_rA^{-p} \notag \\
& S^3(p,q,r)=S_pS_qS_r \notag \\
& S^0(p,q,r)=A^{-p-q-r} \notag
\end{align}

\begin{lm} The Kauffman bracket for $\K$ is
$$
\begin{aligned}
&\langle\K\rangle
=\left(S^0(p,q,r)S^0(a,b,c) + S^2(p,q,r)S^2(a,b,c)\right.\\
&\quad\qquad\qquad\qquad+S_pS_aA^{-q-r-b-c} + S_pS_cA^{-q-r-a-b}+ S_qS_aA^{-p-r-b-c}\\
&\quad\qquad\qquad\qquad+\left.S_qS_bA^{-p-r-a-c}+S_rS_bA^{-p-q-a-c}
+
S_rS_cA^{-p-q-a-b}\right)\\
&+\left(S^1(p,q,r)S^0(a,b,c) +S^0(p,q,r)S^1(a,b,c) + S^2(p,q,r)S^1(a,b,c)
\right.\\
&\left.\qquad +S^1(p,q,r)S^2(a,b,c)+ S^3(p,q,r)S^2(a,b,c)+
S^2(p,q,r)S^3(a,b,c)\right)(-A^{-2} -A^2)\\
&+\left(S^2(p,q,r)S^0(a,b,c) + S^0(p,q,r)S^2(a,b,c) +
S^3(p,q,r)S^1(a,b,c)\right.\\
&\qquad+ S^1(p,q,r)S^3(a,b,c) +S^3(p,q,r)S^3(a,b,c) \\
&\left.\qquad+S_pS_bA^{-q-r-a-c} +S_qS_cA^{-p-r-a-b}
+S_rS_aA^{-p-q-b-c}\right)(-A^{-2}-A^2)^2 \\
&+\left(S^3(p,q,r)S^0(a,b,c) + S^0(p,q,r)S^3(a,b,c)\right)(-A^{-2}-A^2)^3
\notag
\end{aligned}
$$
\end{lm}

\begin{proof} This is by a straightforward computation (with great patience).
\end{proof}

So we see that the Kauffman bracket for $\K$ is mostly symmetric with
respect to
permutations of $\{p,q,r\}$ and $\{a,b,c\}$ except for the
\begin{align}
& S_pS_aA^{-q-r-b-c} + S_pS_cA^{-q-r-a-b} + S_qS_aA^{-p-r-b-c}
\notag \\ &+S_qS_bA^{-p-r-a-c} +S_rS_bA^{-p-q-a-c}
+S_rS_cA^{-p-q-a-b} \notag \end{align}
 and $$\left(S_pS_bA^{-q-r-a-c}
+S_qS_cA^{-p-r-a-b} +S_rS_aA^{-p-q-b-c}\right)(-A^{-2}-A^2)^2$$ terms.

\begin{lm} Let $w=p+q+r+a+b+c$. We have:
\begin{equation}
\langle\K\rangle-\langle K\left(\begin{smallmatrix} p & q & r\\
b & a & c
\end{smallmatrix}\right)\rangle=A^{-w}(S_pA^p-S_rA^r)(S_aA^a-S_bA^b)\,(1-(-A^2-A^{-1})^2)
\end{equation}

\begin{equation}
\langle\K\rangle-\langle K\left(\begin{smallmatrix} p & q & r\\
a & c & b
\end{smallmatrix}\right)\rangle=A^{-w}(S_pA^p-S_qA^q)(S_cA^c-S_bA^b)\,(1-(-A^2-A^{-1})^2)
\end{equation}

\begin{equation}
\langle\K\rangle-\langle K\left(\begin{smallmatrix} p & q & r\\
c & b & a
\end{smallmatrix}\right)\rangle=A^{-w}(S_qA^q-S_rA^r)(S_aA^a-S_cA^c)\,(1-(-A^2-A^{-1})^2)
\end{equation}

\begin{equation}
\langle\K\rangle-\langle K\left(\begin{smallmatrix} p & q & r\\
c & a & b
\end{smallmatrix}\right)\rangle=A^{-w}
\left|\begin{matrix} S_pA^p & S_qA^q & S_rA^r \\
                     S_cA^c & S_aA^a & S_bA^b \\
                     1 & 1 & 1 \end{matrix}\right|\,(1-(-A^2-A^{-1})^2)
\end{equation}

\begin{equation}\langle\K\rangle-\langle K\left(\begin{smallmatrix} p & q & r\\
b & c & a
\end{smallmatrix}\right)\rangle=-A^{-w}
\left|\begin{matrix} S_pA^p & S_qA^q & S_rA^r \\
                     S_aA^a & S_bA^b & S_cA^c \\
                     1 & 1 & 1 \end{matrix}\right|\,(1-(-A^2-A^{-1})^2)
\end{equation}
\end{lm}

Consider now the special case that all $p,q,r,a,b,c$ are even.
In this case, the writhe of the knot $\K$ is $w=p+q+r+a+b+c$.
Notice that $w$ is also the writhe of $K\left(\begin{smallmatrix} p & q & r\\
c & a & b
\end{smallmatrix}\right)$ as well as the writhe of
$K\left(\begin{smallmatrix} p & q & r\\
b & c & a\end{smallmatrix}\right)$. So we have the following results.

\begin{lm} Let $\tau\in S_3$ be a transposition. For $p,q,r,a,b,c$
even,
$$J_{\K}=J_{K\left(\begin{smallmatrix} p & q & r\\
\tau(a) & \tau(b) & \tau(c)
\end{smallmatrix}\right)}$$
iff
$$\K=K\left(\begin{smallmatrix} p & q & r\\
\tau(a) & \tau(b) & \tau(c)
\end{smallmatrix}\right)$$
by a $D_3$ symmetry.
\end{lm}

\begin{thm} Suppose $p,q,r,a,b,c$ are even integers.
If $\K=K\left(\begin{smallmatrix} p & q & r\\
c & a & b
\end{smallmatrix}\right)$, then
$$\left|\begin{matrix} S_pA^p & S_qA^q & S_rA^r \\
                     S_cA^c & S_aA^a & S_bA^b \\
                     1 & 1 & 1 \end{matrix}\right|=0.
                    $$

Also, if $\K=K\left(\begin{smallmatrix} p & q & r\\
b & c & a\end{smallmatrix}\right)$, then
$$\left|\begin{matrix} S_pA^p & S_qA^q & S_rA^r \\
                     S_aA^a & S_bA^b & S_cA^c \\
                     1 & 1 & 1 \end{matrix}\right|=0.$$
\end{thm}

\begin{example}{\rm Consider $K\left(\begin{smallmatrix} 4 & 8 & 12\\
4 & 6 & 2
\end{smallmatrix}\right)$ and $K\left(\begin{smallmatrix} 4 & 8 & 12\\
2 & 4 & 6
\end{smallmatrix}\right)$. Since
$$\left|\begin{matrix} 4 & 8 & 12 \\
                       2& 4 & 6\\ 1 & 1 & 1 \end{matrix}\right|=0,$$
we can not use the Conway polynomial to distinguish these two knots.
But
$$\left|\begin{matrix} S_4A^4 & S_8A^8 & S_{12}A^{12}\\
                       S_2A^2 & S_4A^4 & S_6A^6 \\ 1& 1& 1 \end{matrix}\right|
=\frac{A^{32}-2A^{40}+2A^{56}-A^{64}}{(A^2+A^{-2})^2}\neq 0.$$
Thus,
$$K\left(\begin{smallmatrix} 4 & 8 & 12\\
4 & 6 & 2
\end{smallmatrix}\right)\neq K\left(\begin{smallmatrix} 4 & 8 & 12\\
2 & 4 & 6
\end{smallmatrix}\right).$$}
\end{example}

\subsection{Other permutations of $p,q,r,a,b,c$} Use the same method, We consider briefly one case of a permutation of $p,q,r,a,b,c$ and its effect on the knot $\K$.

Due to the limitation of our method, we consider only the case
that all $p,q,r,a,b,c$ are even so that the writhe of $\K$
is $w=p+q+r+a+b+c$ and it will not change if we permute $p,q,r,a,b,c$.
The question is, in this case, what happens if
$\langle \K \rangle = \langle K\left(\begin{smallmatrix} a & q &r\\
p & b & c
\end{smallmatrix}\right) \rangle $?  The difference of the Kauffman bracket
is:
$$\begin{aligned}
\langle K\left(\begin{smallmatrix} p & q &r\\
a & b & c
\end{smallmatrix}\right) \rangle -\langle K\left(\begin{smallmatrix} a & q &r\\
p & b & c
\end{smallmatrix}\right) \rangle
=&(S_pA^{-q-r-a-b-c}-S_aA^{-p-q-r-b-c})\\
 &\cdot(S_bS_cA^{-p-q-r-a}-S_qS_rA^{-p-a-b-c})\\
&\cdot(-A^{-2} - A^2 - (-A^{-2}-A^2)^3) \\
&+(S_qS_bS_cA^{-p-r-a}+S_rS_bS_cA^{-p-q-a} + S_cA^{-p-q-r-a-b}\\
&\quad- S_qS_rS_bA^{-p-a-c} -S_qS_rS_cA^{-p-a-b} -
S_qA^{-p-r-a-b-c})\\
&\cdot(1-(-A^{-2}-A^2)^2).
\end{aligned}
$$

If $p=a$, then we certainly have $\K=K\left(\begin{smallmatrix} a & q &r\\
p & b & c
\end{smallmatrix}\right)$.  So consider the case when
$p\neq a$.  Then, if $\langle \K \rangle = \langle K\left(\begin{smallmatrix} a & q &r\\
p & b & c
\end{smallmatrix}\right) \rangle $, we have
\begin{equation}\label{iden}
\begin{aligned}
&S_qS_bS_cA^{-r}+S_rS_bS_cA^{-q}+S_cA^{-q-r-b}-S_qS_rS_bA^{-c}-
S_qS_rS_cA^{-b}-S_qA^{-r-b-c}\\
&=(A^2+A^{-2})(S_bS_cA^{-q-r}-S_qS_rA^{-b-c}).
\end{aligned}
\end{equation}

We have
$$\begin{aligned}
S_q&=A^{2-q}-A^{2-q+4}+A^{2-q+8}-\cdots+(-1)^{q-1}A^{3q-2}\\
&=A^{-q}(A^2-A^6+A^{10}+\dots+(-1)^{q-1}A^{2+4(q-1)}).
\end{aligned}$$
Denote
$$\hat{S}_q:=S_qA^q=A^2-A^6+A^{10}+\dots+(-1)^{q-1}A^{2+4(q-1)}=
\frac{1-A^{4q}}{A^2+A^{-2}}.$$

We can rewrite Equation (\ref{iden}) as follows:

\begin{equation}\label{iden2}\begin{aligned}
&\hat{S}_q\hat{S}_b\hat{S}_c+\hat{S}_r\hat{S}_b\hat{S}_c-
\hat{S}_q\hat{S}_r\hat{S}_b-
\hat{S}_q\hat{S}_r\hat{S}_c\\
&=(A^2+A^{-2})(\hat{S}_b\hat{S}_c-\hat{S}_q\hat{S}_r)+\hat{S}_q-
\hat{S}_c
\end{aligned}
\end{equation}

Put in the fraction form of $\hat{S}_q$, we get
\begin{equation}\label{iden3}
\begin{aligned}
&(1-A^{4q})(1-A^{4b})(1-A^{4c})+(1-A^{4r})(1-A^{4b})(1-A^{4c})\\
      &-(1-A^{4q})(1-A^{4r})(1-A^{4b})-(1-A^{4q})(1-A^{4r})(1-A^{4c})\\
&=(A^2+A^{-2})^2\left[(1-A^{4b})(1-A^{4c})-(1-A^{4q})(1-A^{4r})
+A^{4c}-A^{4q}\right].
\end{aligned}
\end{equation}

Let us take derivative with respect to $A$ of the both sides of Equation (\ref{iden3}) at $A=1$.
The left hand side is zero, and the right hand side is zero only
when $q=c$.

So now we assume $q=c$. Then the equation (3) becomes
$$(1-A^{4q})^2(A^{4b}-A^{4r})=(A^2+A^{-2})^2(1-A^{4q})(A^{4b}-A^{4r}).$$
This is true only when $q=0$ or $b=r$. So we have proved the
following theorem.

\begin{thm} Suppose that all $p,q,r,a,b,c$ are even integers and $p\neq a$.
If
$$\K=K\left(\begin{smallmatrix} a & q &r\\
p & b & c
\end{smallmatrix}\right),$$
then either (1) $q=c=0$, or (2) $q=c$ and $b=r$.
\end{thm}

We finish this section by proposing two problems.

\begin{prob} Classify knots $\K$ for $p,q,r,a,b,c$ all even.
\end{prob}

\begin{prob} Show that the knot $8_{18}$ has girth $g>3$.
\end{prob}

\section{Planar Tree Pair Representations of Knots}

The following table lists planar tree pair
representations of knots with girth
$g\leq 3$ from Rolfsen's knots and links table \cite{R}. The
representations are not unique.
\smallskip

\begin{tabular}{c l}
$3_1$ & (3) \\
$4_1$ & (2,-2) \\
$5_1$ & (5) \\
$5_2$ & (2,-3) \\
$6_1$ & (2,-4) \\
$6_2$ & $\left(\begin{smallmatrix} 0 & 2 & -2\\ 0 & -1 & -1
\end{smallmatrix}
\right)$ \\
$6_3$ & $\left(\begin{smallmatrix} 2 & 0 & 1\\ -1 & -1 & -1
\end{smallmatrix}
\right)$ \\
$7_1$ & (7) \\
$7_2$ & (2,-5) \\
$7_3$ & (3,-4) \\
$7_4$ & $\left(\begin{smallmatrix} -1 & -1 & 0\\ 1 & 2 & 2
\end{smallmatrix}
\right)$ \\
$7_5$ & $\left(\begin{smallmatrix} 3 & 1 & 0\\ -1 & -0 & -2
\end{smallmatrix}
\right)$ \\
$7_6$ & $\left(\begin{smallmatrix} -1 & 2 & 0\\ 1 & 2 & 2
\end{smallmatrix}
\right)$ \\
$7_7$ & $\left(\begin{smallmatrix} 2 & 2 & 0\\ -1 & -1 & -1
\end{smallmatrix}
\right)$ \\
$8_1$ & (2,-6) \\
$8_2$ & $\left(\begin{smallmatrix} 0 & -1 & -1\\ 1 & 1 & 4
\end{smallmatrix}
\right)$ \\
$8_3$ & (4,-4) \\
$8_4$ & $\left(\begin{smallmatrix} -1 & -1 & 0\\ 1 & 2 & 3
\end{smallmatrix}
\right)$ \\
$8_5$ & $\left(\begin{smallmatrix} 0 & 1 & 1\\ -1 & -3 & -2
\end{smallmatrix}
\right)$ \\
$8_6$ & $\left(\begin{smallmatrix} 1 & 3 & 0\\ -1 & -3 & 0
\end{smallmatrix}
\right)$ \\
$8_7$ & $\left(\begin{smallmatrix} 0 & -1 & -4\\ 1 & 1 & 1
\end{smallmatrix}
\right)$ \\
$8_8$ & $\left(\begin{smallmatrix} 0 & -2 & -1\\ 3 & 1 & 1
\end{smallmatrix}
\right)$ \\
$8_9$ & $\left(\begin{smallmatrix} 3 & 0 & 1\\ -1 & -2 & -1
\end{smallmatrix}
\right)$ \\
$8_{10}$ & $\left(\begin{smallmatrix} 1 & 0 & 2\\ -1 & -1 & -3
\end{smallmatrix}
\right)$ \\
$8_{11}$ & $\left(\begin{smallmatrix} -3 & -1 & 0\\ 1 & 1 & 2
\end{smallmatrix}
\right)$
\end{tabular}
\begin{tabular}{c l}

$8_{12}$ & $\left(\begin{smallmatrix} -2 & -2 & 0\\ 2 & 2 & 0
\end{smallmatrix}
\right)$ \\
$8_{13}$ & $\left(\begin{smallmatrix} 2 & 0 & 3\\ -1 & -1 & -1
\end{smallmatrix}
\right)$ \\
$8_{14}$ & $\left(\begin{smallmatrix} 2 & 2 & 0\\ -1 & -1 & -2
\end{smallmatrix}
\right)$ \\
$8_{15}$ & $\left(\begin{smallmatrix} 2 & 2 & 0\\ -2 & -1 & -1
\end{smallmatrix}
\right)$ \\
$8_{16}$ & $\left(\begin{smallmatrix} 2 & 2 & 1\\ -1 & -1 & -1
\end{smallmatrix}
\right)$ \\
$8_{17}$ & $\left(\begin{smallmatrix} -2 & -1 & -1\\ 1 & 1 & 2
\end{smallmatrix}
\right)$ \\
$8_{18}$ & ? \\
$8_{19}$ & $\left(\begin{smallmatrix} -2 & -1 & -1\\ 1 & -1 & -2
\end{smallmatrix}
\right)$ \\
$8_{20}$ & $\left(\begin{smallmatrix} 2 & -1 & -1\\ -1 & 1 & -2
\end{smallmatrix}
\right)$ \\
$8_{21}$ & $\left(\begin{smallmatrix} 2 & 1 & 1\\ 1 & 1 & 2
\end{smallmatrix}
\right)$ \\

$9_1$ & (9) \\
$9_2$ & (2,-7) \\
$9_3$ & (3,-6) \\
$9_4$ & (4,-5) \\
$9_5$ & $\left(\begin{smallmatrix} -1 & -1 & 0\\ 1 & 2 & 4
\end{smallmatrix}
\right)$ \\

$9_6$ & $\left(\begin{smallmatrix} 5 & 1 & 0\\ -1 & 0 & -2
\end{smallmatrix}
\right)$ \\
$9_7$ & $\left(\begin{smallmatrix} 3 & 1 & 0\\ -1 & 0 & -4
\end{smallmatrix}
\right)$ \\
$9_8$ & $\left(\begin{smallmatrix} 0 & -1 & -2\\ 1 & 1 & 4
\end{smallmatrix}
\right)$ \\
$9_9$ & $\left(\begin{smallmatrix} -1 & -2 & 0\\ 3 & 0 & 3
\end{smallmatrix}
\right)$ \\
$9_{10}$ & $\left(\begin{smallmatrix} -1 & -3 & 0\\ 3 & 0 & 2
\end{smallmatrix}
\right)$ \\
$9_{11}$ & $\left(\begin{smallmatrix} 1 & 2 & 0\\ -1 & -2 & -3
\end{smallmatrix}
\right)$ \\
$9_{12}$ & $\left(\begin{smallmatrix} -1 & -4 & 0\\ 1 & 2 & 1
\end{smallmatrix}
\right)$ \\
$9_{13}$ & $\left(\begin{smallmatrix} -1 & -3 & 0\\ 1 & 2 & 2
\end{smallmatrix}
\right)$ \\
$9_{14}$ & $\left(\begin{smallmatrix} 0 & 2 & 4\\ -1 & -1 & -1
\end{smallmatrix}
\right)$ \\
$9_{18}$ & $\left(\begin{smallmatrix} -3 & 0 & -2\\ 2 & 0 & 2
\end{smallmatrix}
\right)$
\end{tabular}
\begin{tabular}{c l}
$9_{19}$ & $\left(\begin{smallmatrix} 2 & 2 & 0\\
-1 & -1 & -3
\end{smallmatrix}
\right)$ \\
$9_{21}$ & $\left(\begin{smallmatrix} 2 & 3 & 0\\ -1 & -1 & -2
\end{smallmatrix}
\right)$ \\
$9_{23}$ &  $\left(\begin{smallmatrix} 0 & 2 & 2\\ -2 & -1 & -2
\end{smallmatrix}
\right)$\\
$9_{25}$ & $\left(\begin{smallmatrix} 0 & 2 & 2\\ -2 & -2 & -1
\end{smallmatrix}
\right)$\\
$9_{35}$ & $\left(\begin{smallmatrix} 1 & 1 & 0\\ -3 & -2 & -2
\end{smallmatrix}
\right)$ \\
$9_{36}$ & $\left(\begin{smallmatrix} 1 & 2 & 0\\ -2 & -2 & -2
\end{smallmatrix}
\right)$ \\
$9_{37}$ & $\left(\begin{smallmatrix} 2 & 2 & 0\\ -3 & -1 & -1
\end{smallmatrix}
\right)$ \\
$9_{39}$ & $\left(\begin{smallmatrix} 2 & 2 & 1\\ -2 & -1 & -1
\end{smallmatrix}
\right)$\\
$9_{41}$ & $\left(\begin{smallmatrix} -2 & -2 & -2\\ 1 & 1 & 1
\end{smallmatrix}
\right)$ \\
$9_{42}$ & $\left(\begin{smallmatrix} 0 & 1 & 2\\ -2 & 2 & -2
\end{smallmatrix}
\right)$ \\
$9_{46}$ & $\left(\begin{smallmatrix} -1 & 1 & 0\\ -3 & -2 & 2
\end{smallmatrix}
\right)$ \\
$9_{48}$ & $\left(\begin{smallmatrix} 2 & 2 & 1\\ 2 & 1 & 1
\end{smallmatrix}
\right)$ \\
$9_{49}$ & $\left(\begin{smallmatrix} 2 & 2 & 1\\ 2 & -1 & -1
\end{smallmatrix}
\right)$ \\
$10_1$ & (2,-8) \\
$10_2$ & $\left(\begin{smallmatrix} -1 & 0 & -1\\ 0 & 6 & 2
\end{smallmatrix}
\right)$ \\
$10_3$ & (4,-6) \\
$10_4$ & $\left(\begin{smallmatrix} 0 & -1 & -1\\ 1 & 2 & 5
\end{smallmatrix}
\right)$ \\
$10_5$ & $\left(\begin{smallmatrix} 2 & 0 & 1\\ -1 & -5 & -1
\end{smallmatrix}
\right)$ \\
$10_6$ & $\left(\begin{smallmatrix} 5 & 1 & 0\\ -1 & 0 & -3
\end{smallmatrix}
\right)$ \\
$10_7$ & $\left(\begin{smallmatrix} -5 & -1 & 0\\ 2 & 0 & 2
\end{smallmatrix}
\right)$\\
$10_8$ & $\left(\begin{smallmatrix} -1 & -1 & 0\\ 1 & 3 & 4
\end{smallmatrix}
\right)$ \\
$10_9$ & $\left(\begin{smallmatrix} 1 & 3 & 0\\ -1 & -1 & -4
\end{smallmatrix}
\right)$ \\
$10_{11}$ & $\left(\begin{smallmatrix} 0 & 3 & 3\\ -3 & -1 & 0
\end{smallmatrix}
\right)$ \\
$10_{12}$ & $\left(\begin{smallmatrix} -2 & -1 & -2\\ 1 & 1 & 3
\end{smallmatrix}
\right)$ \\
$10_{13}$ & $\left(\begin{smallmatrix} 0 & 2 & 2\\ -2 & -4 & 0
\end{smallmatrix}
\right)$
\end{tabular}
\begin{tabular}{c l}
$10_{14}$ & $\left(\begin{smallmatrix} 4 & 2 & 0\\ -1 & -1 & -2
\end{smallmatrix}
\right)$ \\
$10_{15}$ & $\left(\begin{smallmatrix} 2 & 0 & 1\\ -3 & -3 & -1
\end{smallmatrix}
\right)$ \\
$10_{16}$ & $\left(\begin{smallmatrix} 3 & 0 & 1\\ -2 & -3 & -1
\end{smallmatrix}
\right)$ \\
$10_{17}$ & $\left(\begin{smallmatrix} 4 & 0 & 1\\ -1 & -3 & 1
\end{smallmatrix}
\right)$ \\
$10_{18}$ & $\left(\begin{smallmatrix} 2 & 0 & 4\\ -2 & -1 & -1
\end{smallmatrix}
\right)$ \\
$10_{20}$ & $\left(\begin{smallmatrix} 1 & 3 & 0\\ -1 & -5 & 0
\end{smallmatrix}
\right)$ \\
$10_{21}$ & $\left(\begin{smallmatrix} 0 & -1 & -3\\ 1 & 1 & 4
\end{smallmatrix}
\right)$ \\
$10_{22}$ & $\left(\begin{smallmatrix} -3 & -1 & 0\\ 3 & 0 & 3
\end{smallmatrix}
\right)$ \\
$10_{23}$ & $\left(\begin{smallmatrix} 2 & 0 & 3\\ -1 & -3 & -1
\end{smallmatrix}
\right)$ \\
$10_{24}$ & $\left(\begin{smallmatrix} 0 & 2 & 2\\ -3 & -3 & 0
\end{smallmatrix}
\right)$ \\
$10_{31}$ & $\left(\begin{smallmatrix} 3 & 0 & 2\\ -1 & -3 & -1
\end{smallmatrix}
\right)$ \\
$10_{34}$ & $\left(\begin{smallmatrix} -2 & 0 & -1\\ 5 & 1 & 1
\end{smallmatrix}
\right)$ \\
$10_{35}$ & $\left(\begin{smallmatrix} 0 & -2 & -2\\ 4 & 2 & 0
\end{smallmatrix}
\right)$ \\
$10_{36}$ & $\left(\begin{smallmatrix} 0 & 2 & 2\\ -4 & -1 & -1
\end{smallmatrix}
\right)$ \\
$10_{37}$ & $\left(\begin{smallmatrix} 3 & 0 & 2\\ 0 & -3 & -2
\end{smallmatrix}
\right)$ \\
$10_{46}$ & $\left(\begin{smallmatrix} 0 & 1 & 1\\ -1 & -3 & -4
\end{smallmatrix}
\right)$ \\
$10_{47}$ & $\left(\begin{smallmatrix} 1 & 2 & 0\\ -2 & -1 & -4
\end{smallmatrix}
\right)$ \\
$10_{48}$ & $\left(\begin{smallmatrix} 1 & 0 & 4\\ -1 & -1 & -3
\end{smallmatrix}
\right)$ \\
$10_{50}$ & $\left(\begin{smallmatrix} 1 & 0 & 3\\ -1 & -2 & -3
\end{smallmatrix}
\right)$ \\
$10_{54}$ & $\left(\begin{smallmatrix} 1 & 0 & 2\\ -1 & -3 & -3
\end{smallmatrix}
\right)$ \\
$10_{55}$ & $\left(\begin{smallmatrix} 2 & 2 & 0\\ -2 & -3 & -1
\end{smallmatrix}
\right)$ \\
$10_{58}$ & $\left(\begin{smallmatrix} 0 & 2 & 2\\ -2 & -2 & -2
\end{smallmatrix}
\right)$ \\
$10_{61}$ & $\left(\begin{smallmatrix} 1 & 1 & 0\\ -3 & -2 & -3
\end{smallmatrix}
\right)$ \\
$10_{62}$ & $\left(\begin{smallmatrix} 1 & 2 & 0\\ -3 & -1 & -3
\end{smallmatrix}
\right)$ \\
$10_{63}$ & $\left(\begin{smallmatrix} 2 & 0 & 2\\ -1 & -1 & -4
\end{smallmatrix}
\right)$
\end{tabular}

\begin{tabular}{c l}
$10_{64}$ & $\left(\begin{smallmatrix} 0 & 3 & 1\\ -1 & -3 & -2
\end{smallmatrix}
\right)$ \\
$10_{67}$ & $\left(\begin{smallmatrix} 2 & 2 & 0\\ -3 & -1 & -2
\end{smallmatrix}
\right)$ \\
$10_{82}$ & $\left(\begin{smallmatrix} 4 & 1 & 1\\ -1 & -1 & -2
\end{smallmatrix}
\right)$ \\
$10_{85}$ & $\left(\begin{smallmatrix} 1 & 4 & 2\\ -1 & -1 & -1
\end{smallmatrix}
\right)$ \\
$10_{91}$ & $\left(\begin{smallmatrix} 3 & 1 & 1\\ -1 & -2 & -2
\end{smallmatrix}
\right)$ \\
$10_{94}$ & $\left(\begin{smallmatrix} 2 & 1 & 1\\ -1 & -3 & -2
\end{smallmatrix}
\right)$ \\
$10_{99}$ & $\left(\begin{smallmatrix} 2 & 2 & 1\\ -1 & -2 & -2
\end{smallmatrix}
\right)$ \\
$10_{102}$ & $\left(\begin{smallmatrix} -2 & -1 & -1\\ -2 & -1 & -3
\end{smallmatrix}
\right)$ \\
$10_{108}$ & $\left(\begin{smallmatrix} -2 & -3 & -2\\ 1 & 1 & 1
\end{smallmatrix}
\right)$ \\
$10_{124}$ & $\left(\begin{smallmatrix} 0 & -1 & 1\\ 1 & -3 & -4
\end{smallmatrix}
\right)$ \\
$10_{125}$ & $\left(\begin{smallmatrix} 0 & 1 & -1\\ -1 & 3 & -4
\end{smallmatrix}
\right)$ \\
$10_{126}$ & $\left(\begin{smallmatrix} 0 & 1 & -1\\ -1 & -3 & 4
\end{smallmatrix}
\right)$ \\
$10_{128}$ & $\left(\begin{smallmatrix} -2 & -2 & 1\\ 3 & -1 & -1
\end{smallmatrix}
\right)$ \\
$10_{129}$ & $\left(\begin{smallmatrix} 2 & -2 & -1\\ 3 & 1 & 1
\end{smallmatrix}
\right)$ \\
$10_{140}$ & $\left(\begin{smallmatrix} 0 & -1 & 1\\ 2 & -3 & -3
\end{smallmatrix}
\right)$ \\
$10_{142}$ & $\left(\begin{smallmatrix} 0 & -1 & 1\\ 3 & -3 & -2
\end{smallmatrix}
\right)$ \\
$10_{143}$ & $\left(\begin{smallmatrix} 2 & -2 & -1\\ 1 & 2 & 2
\end{smallmatrix}
\right)$ \\
$10_{157}$ & $\left(\begin{smallmatrix} -2 & -2 & -1\\ -1 & -2 & -2
\end{smallmatrix}
\right)$ \\
$10_{158}$ & $\left(\begin{smallmatrix} 2 & 2 & 1\\ 3 & -1 & -1
\end{smallmatrix}
\right)$ \\
$10_{163}$ & $\left(\begin{smallmatrix} 1 & 1 & 1\\ 2 & 2 & 3
\end{smallmatrix}
\right)$
\end{tabular}
\begin{tabular}{c l}
$2_1^2$ & (2) \\
$4_1^2$ & (4) \\
$5_1^2$ & $\left(\begin{smallmatrix} 0 & -1 & -1\\ 1 & 1 & 1
\end{smallmatrix}
\right)$ \\
$6_1^2$ & (6) \\
$6_2^2$ & (3,-3) \\
$6_3^2$ & $\left(\begin{smallmatrix} 1 & 2 & 0\\ -1 & -2 & 0
\end{smallmatrix}
\right)$ \\
$7_1^2$ & $\left(\begin{smallmatrix} -1 & 0 & -1\\ 1 & 0 & 4
\end{smallmatrix}
\right)$ \\
$7_2^2$ & $\left(\begin{smallmatrix} -3 & -1 & 0\\ 2 & 0 & 1
\end{smallmatrix}
\right)$ \\
$7_3^2$ & $\left(\begin{smallmatrix} 1 & 2 & 0\\ -1 & -3 & 0
\end{smallmatrix}
\right)$\\

$7_4^2$ & $\left(\begin{smallmatrix} 0 & 1 & 1\\ -1 & -2 & -2
\end{smallmatrix}
\right)$ \\
$7_5^2$ & $\left(\begin{smallmatrix} 0 & 2 & 1\\ -1 & -2 & -1
\end{smallmatrix}
\right)$ \\
$7_6^2$ & $\left(\begin{smallmatrix} 1 & 2 & 1\\ -1 & -1 & -1
\end{smallmatrix}
\right)$ \\
$7_7^2$ & $\left(\begin{smallmatrix} 0 & 1 & -1\\ 1 & -2 & -2
\end{smallmatrix}
\right)$ \\
$7_8^2$ & $\left(\begin{smallmatrix} 0 & 2 & -1\\ -1 & 2 & -1
\end{smallmatrix}
\right)$ \\
$8_1^2$ & (8) \\
$8_2^2$ & (3,-5) \\
$8_3^2$ & $\left(\begin{smallmatrix} 0 & 2 & 3\\ -2 & -1 & 0
\end{smallmatrix}
\right)$\\
$8_4^2$ & $\left(\begin{smallmatrix} -1 & -2 & 0\\ 3 & 0 & 2
\end{smallmatrix}
\right)$ \\
$8_6^2$ & $\left(\begin{smallmatrix} 1 & 1 & 0\\ -4 & -1 & -1
\end{smallmatrix}
\right)$ \\
$8_9^2$ & $\left(\begin{smallmatrix} 1 & 1 & 0\\ -2 & -2 & -2
\end{smallmatrix}
\right)$
\end{tabular}
\begin{tabular}{c l}
$8_{15}^2$ & $\left(\begin{smallmatrix} -1 & 0 & 2\\ 1 & -2 & -2
\end{smallmatrix}
\right)$ \\
$9_1^2$ & $\left(\begin{smallmatrix} 0 & -1 & -1\\ 1 & 1 & 5
\end{smallmatrix}
\right)$ \\
$9_2^2$ & $\left(\begin{smallmatrix} 2 & 0 & 1\\ -1 & -4 & -1
\end{smallmatrix}
\right)$ \\
$9_3^2$ & $\left(\begin{smallmatrix} 1 & 4 & 0\\ -1 & -3 & 0
\end{smallmatrix}
\right)$ \\
$9_4^2$ & $\left(\begin{smallmatrix} 0 & -1 & -1\\ 3 & 1 & 3
\end{smallmatrix}
\right)$ \\
$9_5^2$ & $\left(\begin{smallmatrix} 3 & 0 & 1\\ -1 & -3 & -1
\end{smallmatrix}
\right)$ \\
$9_6^2$ & $\left(\begin{smallmatrix} 0 & -1 & -3\\ 1 & 1 & 3
\end{smallmatrix}
\right)$ \\
$9_7^2$ & $\left(\begin{smallmatrix} 0 & 2 & 3\\ -1 & -1 & -2
\end{smallmatrix}
\right)$ \\
$9_9^2$ & $\left(\begin{smallmatrix} 3 & 3 & 0\\ -1 & -1 & -1
\end{smallmatrix}
\right)$ \\
$9_{10}^2$ & $\left(\begin{smallmatrix} 1 & 2 & 0\\ -1 & -5 & 0
\end{smallmatrix}
\right)$ \\
$9_{13}^2$ & $\left(\begin{smallmatrix} 0 & 1 & 1\\ -1 & -2 & -4
\end{smallmatrix}
\right)$ \\
$9_{14}^2$ & $\left(\begin{smallmatrix} 0 & 4 & 1\\ -1 & -2 & -1
\end{smallmatrix}
\right)$ \\
$9_{15}^2$ & $\left(\begin{smallmatrix} 0 & 3 & 1\\ -2 & -2 & -1
\end{smallmatrix}
\right)$ \\
$9_{17}^2$ & $\left(\begin{smallmatrix} 0 & 2 & 1\\ -3 & -2 & -1
\end{smallmatrix}
\right)$ \\
$9_{19}^2$ & $\left(\begin{smallmatrix} 0 & 1 & 1\\ -1 & -3 & -3
\end{smallmatrix}
\right)$ \\
$9_{20}^2$ & $\left(\begin{smallmatrix} 0 & 1 & 2\\ -1 & -4 & -1
\end{smallmatrix}
\right)$ \\
$9_{21}^2$ & $\left(\begin{smallmatrix} 1 & 0 & 3\\ -1 & -1 & -3
\end{smallmatrix}
\right)$ \\
$9_{22}^2$ & $\left(\begin{smallmatrix} 2 & 3 & 0\\ -2 & -1 & -1
\end{smallmatrix}
\right)$ \\
$9_{23}^2$ & $\left(\begin{smallmatrix} 0 & 2 & 1\\ -1 & -3 & -2
\end{smallmatrix}
\right)$ \\
$9_{31}^2$ & $\left(\begin{smallmatrix} 4 & 1 & 1\\ -1 & -1 & -1
\end{smallmatrix}
\right)$
\end{tabular}
\begin{tabular}{c l}
$9_{35}^2$ & $\left(\begin{smallmatrix} 1 & 3 & 2\\ -1 & -1 & -1
\end{smallmatrix}
\right)$ \\
$9_{36}^2$ & $\left(\begin{smallmatrix} 1 & 3 & 1\\ -1 & -1 & -2
\end{smallmatrix}
\right)$ \\
$9_{43}^2$ & $\left(\begin{smallmatrix} 0 & 1 & 1\\ -1 & 2 & -4
\end{smallmatrix}
\right)$ \\
$9_{44}^2$ & $\left(\begin{smallmatrix} 4 & 0 & -1\\ -1 & 1 & -2
\end{smallmatrix}
\right)$ \\
$9_{45}^2$ & $\left(\begin{smallmatrix} 3 & 0 & -1\\ -2 & 1 & -2
\end{smallmatrix}
\right)$ \\
$9_{47}^2$ & $\left(\begin{smallmatrix} 2 & 0 & -1\\ -3 & 1 & -2
\end{smallmatrix}
\right)$ \\
$9_{49}^2$ & $\left(\begin{smallmatrix} 0 & -1 & 1\\ 1 & -3 & -3
\end{smallmatrix}
\right)$ \\
$9_{50}^2$ & $\left(\begin{smallmatrix} 2 & 0 & -1\\ -1 & 1 & -4
\end{smallmatrix}
\right)$ \\
$9_{51}^2$ & $\left(\begin{smallmatrix} 3 & -1 & 0\\ -3 & 1 & -1
\end{smallmatrix}
\right)$ \\
$9_{52}^2$ & $\left(\begin{smallmatrix} 2 & 3 & 0\\ 2 & -1 & -1
\end{smallmatrix}
\right)$ \\
$9_{53}^2$ & $\left(\begin{smallmatrix} -2 & 0 & 1\\ 1 & -2 & -3
\end{smallmatrix}
\right)$ \\
$9_{54}^2$ & $\left(\begin{smallmatrix} 2 & 0 & -1\\ -1 & 2 & -3
\end{smallmatrix}
\right)$ \\
$6_2^3$ & $\left(\begin{smallmatrix} -1 & -1 & -1\\ 1 & 1 & 1
\end{smallmatrix}
\right)$ \\
$6_3^3$ & $\left(\begin{smallmatrix} -1 & 0 & 1\\ 1 & -1 & -2
\end{smallmatrix}
\right)$ \\
$8_1^3$ & $\left(\begin{smallmatrix} 0 & 1 & 1\\ -1 & -2 & -3
\end{smallmatrix}
\right)$ \\
$8_2^3$ & $\left(\begin{smallmatrix} 3 & 1 & 0\\ -2 & -1 & -1
\end{smallmatrix}
\right)$ \\
$8_5^3$ & $\left(\begin{smallmatrix} 1 & 1 & 3\\ -1 & -1 & -1
\end{smallmatrix}
\right)$ \\
$8_6^3$ & $\left(\begin{smallmatrix} 1 & 1 & 2\\ -2 & -1 & -1
\end{smallmatrix}
\right)$ \\
$8_7^3$ & $\left(\begin{smallmatrix} 0 & -1 & 1\\ 1 & -2 & -3
\end{smallmatrix}
\right)$ \\
$8_8^3$ & $\left(\begin{smallmatrix} 0 & -1 & 3\\ 1 & -2 & -1
\end{smallmatrix}
\right)$ \\
\end{tabular}


\bigskip





\begin{thebibliography}{BM}

\bibitem{A} C.C. Adams, The Knot Book, W.H. Freeman
(1999).

\bibitem{Conway} J.H. Conway, \textit
{An enumeration of knots and links, and some of their algebraic
properties},
Computational Problems in Abstract Algebra, Pergamon (1970), pp. 329--350.

\bibitem{DowThis1} C.H. Dowker and M.B. Thistlethwaite, \textit{
On the classification of knots},
C. R. Math. Rep. Acad. Sci. Canada {\bf 4}(1982), no. 2, 129--131.

\bibitem{DowThis2} C.H. Dowker and M.B. Thistlethwaite,
\textit{Classification of knot projections},
Topology Appl. {\bf 16}(1983), no. 1, 19--31.

\bibitem{HTW} J. Hoste, M. Thistlethwaite, and J. Weeks,
\textit{The first 1,701,936 knots}, Math. Intelligencer {\bf
20}(1998), no. 4, 33--48.

\bibitem{R} D. Rolfsen, Knots and Links, AMS Chelsea
(2003).




\end{thebibliography}
\end{document}